\documentclass{amsart}
\usepackage{amssymb}
\usepackage{amsmath}
\usepackage{amstext}
\usepackage{amsfonts}
\usepackage{amsthm}
\usepackage{mathrsfs}
\usepackage{centernot}
\usepackage{enumitem}
\usepackage{graphicx}

\newtheorem*{theorem*}{Theorem}
\newtheorem*{cor*}{Corollary}
\newtheorem{theorem}{Theorem}

\begin{document}
\title[Prime geodesic theorem]{On the prime geodesic theorem for hyperbolic $3-$manifolds}
\author{Muharem Avdispahi\'{c}}
\address{University of Sarajevo, Department of Mathematics, Zmaja od Bosne
33-35, 71000 Sarajevo, Bosnia and Herzegovina}
\email{mavdispa@pmf.unsa.ba}

\begin{abstract}
Through the Selberg zeta approach, we reduce the exponent in the error term
of the prime geodesic theorem for cocompact Kleinian groups or Bianchi
groups from Sarnak's $\frac{5}{3}$ to $\frac{3}{2}$. At the cost of
excluding a set of finite logarithmic measure, the bound is further improved
to $\frac{13}{9}$.
\end{abstract}

\maketitle

\section{Introduction}
Let $\mathbb{H}^{3}$ denote the $3-$dimensional hyperbolic space and let $\Gamma $ be
a cofinite subgroup of $PSL\left( 2,%
\mathbb{C}
\right) $. The quotient $\Gamma \diagdown \mathbb{H}^{3}$ is a $3-$dimensional
hyperbolic manifold of finite volume. The prime geodesic theorem in this
setting says that the number $\pi _{\Gamma }\left( x\right) $ of prime
geodesics $P$ with the length $l\left( P\right) \leq \log x$ equals
\begin{equation}\label{e1}
\pi _{\Gamma }\left( x\right) =li\left( x^{2}\right)
+\sum\limits_{n=1}^{M}li\left( x^{s_{n}}\right) +E\left( x\right) \text{,}
\end{equation}
where $s_{1},\ldots ,s_{M}$ are the real zeros of the Selberg zeta function $%
Z_{\Gamma }$ lying in the interval $\left( 1,2\right) $ and $E\left(
x\right) $ is the error term.

For groups of the form $\Gamma =\Gamma _{D}=PSL\left( 2,\mathcal{O}_{K}\right) $,
where $\mathcal{O}_{K}$ is the ring of integers of an imaginary quadratic number field
$K=\mathbb{Q}\left( \sqrt{-D}\right) $ of class number one, Sarnak \cite{S} proved%
\begin{equation}\label{e2}
\pi _{\Gamma }\left( x\right) =li\left( x^{2}\right) +O\left( x^{\frac{5}{3}%
+\varepsilon }\right) \text{.}
\end{equation}

In the particular case of $\Gamma =PSL\left( 2,%
\mathbb{Z}
\left[ i\right] \right) $, Koyama \cite{K1} obtained%
\[
\pi _{\Gamma }\left( x\right) =li\left( x^{2}\right) +O\left( x^{\frac{11}{7}%
+\varepsilon }\right)
\]%
under the mean-Lindel\"{o}f hypothesis.

While Bianchi groups are noncompact, the fact that the contribution of the continuous
spectra is dominated by the contribution of the
discrete spectra enables one to derive \eqref{e2} by the same reasoning that
leads to $E\left( x\right) =O\left( x^{\frac{5}{3}+\varepsilon }\right) $ in
\eqref{e1} in the cocompact case \cite{N1}. The additional
ingredient for achieving \eqref{e2} is the knowledge of the lower bound for the
first eigenvalue of the Laplace-Beltrami operator in the respective setting.

Using the explicit formula for the integrated Chebyshev functions of an
appropriate order, we shall decrease the exponent in the error term $E\left(
x\right) $.

\begin{theorem}\label{th1}
Let $\Gamma \subset PSL\left( 2,%
\mathbb{C}
\right) $ be a cocompact group or a noncompact cofinite group that
satisfies the condition%
\begin{equation}\label{e3}
\sum\limits_{\gamma _{n}>0}\frac{x^{\beta _{n}-1}}{\gamma _{n}^{2}}=O\left(
\frac{1}{1+(\log x)^{3}}\right) \text{ \ }\left( x\rightarrow \infty \right)
\text{,}
\end{equation}
where $\beta _{n}+i\gamma _{n}$ are poles of the scattering determinant.
Then,%
\[
\pi _{\Gamma }\left( x\right) =li\left( x^{2}\right) +\underset{n=1}{\overset%
{M}{\sum }}li\left( x^{s_{n}}\right) +O\left( \frac{x^{\frac{3}{2}}}{\log x}%
\right) \text{ \ }\left( x\rightarrow \infty \right) \text{.}
\]
\end{theorem}

\begin{cor*}
If $\Gamma $ is a Bianchi group, then%
\[
\pi _{\Gamma }\left( x\right) =li\left( x^{2}\right) +O\left( \frac{x^{\frac{%
3}{2}}}{\log x}\right) \text{ \ }\left( x\rightarrow \infty \right) \text{.}
\]
\end{cor*}

The bound in Theorem \ref{th1} and Corollary is the $3-$dimensional analogue of
Randol's $O\left( \frac{x^{\frac{3}{4}}}{\log x}\right) $ in the prime
geodesic theorem for Riemann surfaces \cite{R}. If $\Gamma \subset PSL\left( 2,%
\mathbb{R}
\right) $ is a cocompact Fuchsian group (or, for that matter, a noncompact
cofinite group satisfying an analogue of \eqref{e3} \cite[p. 477]{H}), it is possible to reduce
the exponent in Randol's estimate to $\frac{7}{10}$ outside a set of finite
logarithmic measure \cite{MA1}. Under the generalized Lindel\"{o}f hypothesis, one
can reach $\frac{5}{8}$ in the case of $PSL\left( 2,%
\mathbb{Z}
\right) $, i.e., come half a way between $\frac{3}{4}$ and the expected
exponent $\frac{1}{2}$, outside a set of finite logarithmic measure \cite{MA2}. For a simple proof of the bound $\frac{2}{3}$ in the latter case without Lindel\"{o}f hypothesis, see \cite{MA3}.

Such Gallagherian approach to prime geodesic theorems leads to the following
result in our $3-$dimensional setting.

\begin{theorem}\label{th2}
Let $\Gamma \subset PSL\left( 2,%
\mathbb{C}
\right) $ be a cocompact group or a noncompact congruence group for some
imaginary quadratic number field. Then there exists a set $E$ of finite
logarithmic measure such that
\[
\pi _{\Gamma }\left( x\right) =li\left( x^{2}\right) +\underset{n=1}{\overset%
{M}{\sum }}li\left( x^{s_{n}}\right) +O\left( x^{\frac{13}{9}}\left( \log
x\right) ^{-\frac{7}{9}}\left( \log \log x\right) ^{\frac{2}{9}+\varepsilon
}\right)\text{,}
\]
as $x\rightarrow \infty$, $x\notin E$.
\end{theorem}

\section{Preliminaries}

We use the upper half-space model%
\begin{equation*}
\mathbb{H}^{3}=\left\{ \left( x,y,z\right) \in
\mathbb{R}
^{3}|r>0\right\} =\left\{ \left( z,r\right) |z\in
\mathbb{C}
,r>0\right\} =\left\{ z+rj|r>0\right\}
\end{equation*}
with the hyperbolic metric $ds^{2}=\frac{dx^{2}+dy^{2}+dr^{2}}{r^{2}}$ and
volume form $d\text{v}=\frac{dxdydz}{r^{3}}$. The Laplace-Beltrami operator is
defined by%
\begin{equation*}
\bigtriangleup =-r^{2}\left( \frac{\partial ^{2}}{\partial x^{2}}+\frac{%
\partial ^{2}}{\partial y^{2}}+\frac{\partial ^{2}}{\partial z^{2}}\right) +r%
\frac{\partial }{\partial r}\text{.}
\end{equation*}
The group $PSL\left( 2,%
\mathbb{C}
\right) =SL\left( 2,%
\mathbb{C}
\right) \diagup \left\{ \pm I\right\} $ is the group of orientation
preserving isometries of $\mathbb{H}^{3}$. It acts on $\mathbb{H}^{3}$ transitively by%
\begin{equation*}
\left(
\begin{array}{cc}
a & b \\
c & d%
\end{array}%
\right) \left( v\right) =\frac{\left( az+b\right) \left( \overline{cz+d}%
\right) +a\overline{c}r^{2}+rj}{\left\vert cz+d\right\vert ^{2}+\left\vert
e\right\vert ^{2}r^{2}}\text{,}
\end{equation*}
where we put $v=z+rj$.

Discrete subgroups of $PSL\left( 2,%
\mathbb{C}
\right) $ are known as Kleinian groups. We shall consider cocompact $\Gamma
\subset PSL\left( 2,%
\mathbb{C}
\right) $ as well as a certain class of noncompact cofinite $\Gamma $ in
which cases the quotient space $\Gamma \diagdown \mathbb{H}^{3}$ is a compact
respectively finite volume hyperbolic $3-$manifold.

If the trace $tr\left( P\right) $ of $P\in \Gamma \setminus \left\{
I\right\} $ is real, $P$ is called hyperbolic, parabolic or elliptic
depending whether $\left\vert tr\left( P\right) \right\vert $ is larger,
equal or less than $2$. In all other cases, $P$ is loxodromic.

Every hyperbolic or loxodromic $P\in \Gamma $ is conjugate in $%
PSL\left( 2,%
\mathbb{C}
\right) $ to a unique element%
\begin{gather*}
\begin{pmatrix}
a\left( P\right)  & 0 \\
0 & a\left( P\right) ^{-1}%
\end{pmatrix}%
\text{, }\left\vert a\left( P\right) \right\vert >1\text{.}
\end{gather*}
The norm of $P$ is defined by $N\left( P\right) =\left\vert a\left( P\right)
\right\vert ^{2}$. For $P$ there exist exactly one primitive hyperbolic or
loxodromic $P_{0}\in \Gamma $ and exactly one $n\in
\mathbb{N}
$ such that $P=P_{0}^{n}$.

Based on a correspondence between conjugate classes of $\Gamma $ and free
homotopy classes of closed continuous paths on $\Gamma \diagdown \mathbb{H}%
^{3}$ (see, e.g., \cite{EGM} for necessary details), we are interested in the
number $\pi _{\Gamma }\left( x\right) $ of primitive hyperbolic or
loxodromic conjugate classes $P_{0}$ with the norm $N\left( P_{0}\right)
\leq x$, i.e., we are interested in $\pi _{\Gamma }\left( x\right)
=\sum_{N\left( P_{0}\right) \leq x}1$. This resembles the situation with
the problem of distribution of prime numbers $\pi \left( x\right)
=\sum_{p\leq x}1$ that led Riemann to introduce his famous zeta function.

The Selberg zeta function  $Z_{\Gamma }$ is defined by%
\[
Z_{\Gamma }\left( s\right) =\prod\limits_{P_{0}}\prod\limits_{k,l}\left(
1-a\left( P_{0}\right) ^{-2k}\overline{a\left( P_{0}\right) }^{-2l}N\left(
P_{0}\right) ^{-s}\right) \text{, \ }\text{Re}s>2 \text{,}
\]%
where the first product is over all primitive hyperbolic or loxodromic
conjugacy classes of $\Gamma $ and the second product is over all pairs of
nonnegative integers such that $k\equiv l\left(mod \ m\left(
P_{0}\right) \right) $, $m\left( P\right) $ denoting the order of the
torsion of the centralizer of $P$.

If $\Gamma $ is cocompact, the functional equation for $Z_{\Gamma }$ reads
(see \cite[Cor. 4.4 on p. 209]{EGM})%
\[
Z_{\Gamma }\left( 2-s\right) =\exp \left( -\frac{vol\left( \Gamma \diagdown
\mathbb{H}^{3}\right) }{3\pi }\left( s-1\right) ^{3}+E\left( s-1\right)
\right) Z_{\Gamma }\left( s\right) \text{,}
\]%
where $E=\sum\limits_{R}\frac{\log N\left( P_{0}\right) }{m\left( R\right)
\left\vert tr\left( R\right) ^{2}-4\right\vert }$, the sum being taken over
all elliptic conjugacy classes of $\Gamma $.

In general case, one has (see \cite[Theorem 4.4]{GW})%
\begin{eqnarray*}
Z_{\Gamma }\left( 2-s\right)  &=&Z_{\Gamma }\left( s\right) \left( \frac{%
\Gamma \left( 2-s\right) }{\Gamma \left( s\right) }\right) ^{4\kappa h_{\Gamma }}%
\left[ \varphi \left( 1-s\right) \right] ^{4\kappa}\prod\limits_{k=1}^{l}\left(
\frac{s-1-q_{k}}{1-s-q_{k}}\right) ^{4\kappa b_{k}} \\
&&\exp \left[ \int\limits_{0}^{s-1}4\pi \kappa vol\left( \Gamma \diagdown \mathbb{H%
}^{3}\right) t^{2}dt+\kappa_{1}\left( s-1\right) \right] \text{,}
\end{eqnarray*}%
where $h_{\Gamma}$ is the number of cusps, $\kappa$ and $\kappa_{1}$ are certain constants, $q_{k}$ $\left( 1\leq k\leq
l\right) $ are the poles of the scattering determinant $\varphi $ in $\left(
0,1\right] $ with order $b_{k}$.

(The standard symbol $\Gamma \left( s\right) $ for Euler's gamma function
appearing in the last equation should not cause any confusion with the
notation for the group.)

The relationship between zeros of $Z_{\Gamma }$ and the discrete spectrum of
the Laplace-Beltrami operator on $L^{2}\left( \Gamma \diagdown \mathbb{H}%
^{3}\right) $ is given by $s\left( 2-s\right) =\lambda $. The latter
operator being essentially self-adjoint, one has that $\lambda _{n}\nearrow
+\infty $ $\left( n\rightarrow +\infty \right) $. So, there are finitely
many zeros $s_{0},\ldots ,s_{M}$ lying in $\left( 1,2\right] $ $\left(
\widetilde{s}_{0},\ldots ,\widetilde{s}_{M}\text{ in }\left[ 0,1\right)
\right) $. All others $s_{n}=1\pm it_{n}$, $t_{n}>0$, corresponding to
discrete eigenvalues $\lambda _{n}>1$, are on the critical line $\text{Re}(s)=1
$.

This serves as a ground for the expectation that $E\left( x\right) =O\left(
x^{1+\varepsilon }\right) $ in \eqref{e1}. However, the density of zeros of $%
Z_{\Gamma }$ has prevented all the efforts in establishing the analogue of
von Koch's theorem \cite[p. 84]{I} for Riemann surfaces or higher dimensional
manifolds.

For groups considered in this paper, discrete eigenvalues are distributed
according to
\begin{equation}
N_{\Gamma }\left( T\right) =\sharp \left\{ \lambda _{n}:\lambda _{n}\leq
1+T^{2}\right\} \sim \frac{vol\left( \Gamma \diagdown \mathbb{H}^{3}\right)
}{6\pi ^{2}}T^{3}\text{.}  \label{e34}
\end{equation}%
Namely,  \eqref{e34} is the Weyl law for cocompact $\Gamma $. If $\Gamma $
is noncompact cofinite and satisfies  \eqref{e3}, then $\int\limits_{-T}^{T}%
\frac{\varphi ^{\prime }}{\varphi }\left( 1+it\right) dt=o\left(
T^{2}\right) $, what in combination with the extended Weyl law \cite[Th.
5.4. on p. 307]{EGM}%
\begin{equation*}
N_{\Gamma }\left( T\right) -\frac{1}{4\pi }\int\limits_{-T}^{T}\frac{\varphi
^{\prime }}{\varphi }\left( 1+it\right) dt\sim \frac{vol\left( \Gamma
\diagdown \mathbb{H}^{3}\right) }{6\pi ^{2}}T^{3}
\end{equation*}%
yields \eqref{e34} (cf. \cite[Th. 4.3]{N2}).

An analogue of the classical von Mangoldt function is given by
\begin{equation*}
\Lambda _{\Gamma }\left( P\right) =\frac{N\left( P\right) \log N\left(
P_{0}\right) }{m\left( R\right) \left\vert a\left( P\right) -a\left(
P\right) ^{-1}\right\vert ^{2}}\text{,}
\end{equation*}%
where $P_{0}$ is a primitive element associated to $P$.

Similar to the Riemann zeta case, a convenient tool to study the
distribution of prime geodesics is provided by the Chebyshev function

\begin{equation*}
\psi _{0,\Gamma }\left( x\right) =\sum\limits_{N\left( P\right) \leq
x}\Lambda _{\Gamma }\left( P\right)
\end{equation*}%
and its integrated versions $\psi _{1,\Gamma }\left( x\right)
=\int\limits_{0}^{x}\psi _{0,\Gamma }\left( t\right) dt$, $\psi _{n,\Gamma
}\left( x\right) =\int\limits_{0}^{x}\psi _{n-1,\Gamma }\left( t\right) dt$.

Explicit formulas with an error term for appropriately integrated Chebyshev
function are the starting ground in our proofs.

\section{Proof of Theorem \ref{th1}}\label{s3}
Let $\Gamma $ be a cofinite group satisfying the condition \eqref{e3}. In view of \eqref{e34}, the
explicit formula for $\psi_{3,\Gamma}$ (see \cite[Th. 5.4.]{N2}) can be
written in the form
\begin{gather*}
\psi _{3,\Gamma }\left( x\right) =A_{0}x^{3}+B_{0}x^{3}\log
x+A_{1}x^{2}+A_{2}x+A_{3} \\
+\sum\limits_{n=0}^{M}\frac{x^{s_{n}+3}}{s_{n}\left( s_{n}+1\right) \left(
s_{n}+2\right) \left( s_{n}+3\right) } +\sum\limits_{n=0}^{M}\frac{x^{\widetilde{s}_{n}+3}}{\widetilde{s}_{n}\left( \widetilde{s}_{n}+1\right) \left(\widetilde{s}_{n}+2\right) \left(\widetilde{s}_{n}+3\right) }\\
+\sum\limits_{0\leq t_{n}\leq T}\frac{x^{s_{n}+3}}{s_{n}\left(
s_{n}+1\right) \left( s_{n}+2\right) \left( s_{n}+3\right) }%
+\sum\limits_{0\leq t_{n}\leq T}\frac{x^{\widetilde{s}_{n}+3}}{\widetilde{s}%
_{n}\left( \widetilde{s}_{n}+1\right) \left( \widetilde{s}_{n}+2\right)
\left( \widetilde{s}_{n}+3\right) } \\
+O\left( \frac{x^{4}}{T}\right)
+\sum\limits_{\gamma _{n}\geq 0}\frac{x^{\rho _{n}+3}}{\rho _{n}\left( \rho
_{n}+1\right) \left( \rho _{n}+2\right) \left( \rho _{n}+3\right) }%
 \\ +\sum\limits_{\gamma _{n}\geq 0}\frac{x^{\widetilde{\rho }_{n}+3}}{\widetilde{%
\rho }_{n}\left( \widetilde{\rho }_{n}+1\right) \left( \widetilde{\rho }%
_{n}+2\right) \left( \widetilde{\rho }_{n}+3\right) }\text{,}
\end{gather*}
where $s_{n}=1+it_{n}$, $\widetilde{s}_{n}=1-it_{n}$ are the zeros
of $Z_{\Gamma }$ coming from the discrete spectrum and $\rho _{n}=\beta _{n}+i\gamma
_{n}$, $\widetilde{\rho }_{n}=\beta _{n}-i\gamma_{n}$ are those from the continuous spectrum. (If $\Gamma $ is a cocompact group, then the last two sums in the explicit formula above are obviously void.)

The asymptotics of a nonnegative nondecreasing function $\psi _{0,\Gamma }
$ can be easily derived from the asymptotics of $\psi _{3,\Gamma }$. One
introduces the functions
\begin{gather*}
\Delta _{3}^{+}f\left( x\right)
=\int\limits_{x}^{x+h}\int\limits_{t}^{t+h}\int\limits_{v}^{v+h}f^{\prime
\prime \prime }\left( u\right) dudvdt\\
=f\left( x+3h\right) -3f\left(
x+2h\right) +3f\left( x+h\right) -f\left( x\right)
\end{gather*}
and
\begin{equation*}
\Delta _{3}^{-}f\left( x\right)
=\int\limits_{x-h}^{x}\int\limits_{t-h}^{t}\int\limits_{v-h}^{v}f^{\prime
\prime \prime }\left( u\right) dudvdt\text{.}
\end{equation*}

By the mean value theorem, we get%
\begin{equation*}
h^{-3}\Delta _{3}^{+}\sum\limits_{n=0}^{M}\frac{x^{s_{n}+3}}{s_{n}\left(
s_{n}+1\right) \left( s_{n}+2\right) \left( s_{n}+3\right) }%
=\sum\limits_{n=0}^{M}\frac{x^{s_{n}}}{s_{n}}+O\left( h^{2}\right) \text{.}
\end{equation*}

Further, for the zeros on the critical line, we have
\begin{gather*}
h^{-3}\Delta _{3}^{+}\sum\limits_{\left\vert t_{n}\right\vert \leq T}\frac{%
x^{s_{n}+3}}{s_{n}\left( s_{n}+1\right) \left( s_{n}+2\right) \left(
s_{n}+3\right) }=O\left( x\sum\limits_{\left\vert t_{n}\right\vert \leq T}%
\frac{1}{\left\vert s_{n}\right\vert }\right) =O\left( xT^{2}\right) \text{.}
\end{gather*}

Now,
\begin{gather*}
h^{-3}\Delta _{3}^{+}\sum\limits_{\gamma _{n}\geq 0}\frac{x^{\rho _{n}+3}}{%
\rho _{n}\left( \rho _{n}+1\right) \left( \rho _{n}+2\right) \left( \rho
_{n}+3\right) }\\
=O\left( \frac{1}{h}\sum\limits_{\gamma _{n}\geq 0}\frac{%
x^{\beta _{n}+1}}{\left\vert \rho _{n}\right\vert \left\vert \rho
_{n}+1\right\vert }\right) =O\left( \frac{x^{2}}{h}\right)
\end{gather*}
by \eqref{e3}.

Thus,
\begin{gather*}
\psi _{0,\Gamma }\left( x\right) \leq h^{-3}\Delta _{3}^{+}\psi
_{_{3,\Gamma }}\left( x\right) \\
=\sum\limits_{n=0}^{M}\frac{x^{s_{n}}}{s_{n}}%
+O\left( h^{2}\right) +O\left( xT^{2}\right) +O\left( \frac{x^{4}}{h^{3}T}\right)
+O\left( \frac{x^{2}}{h}\right) \text{.}
\end{gather*}
The optimal choice for the first three $O-$summands on the right-hand side
is $h \sim x^{\frac{3}{4}}$, $T \sim x^{\frac{1}{4}}.$ The fourth term
$O\left( \frac{x^{2}}{h}\right) $ is dominated by the obtained bound $%
O\left( x^{\frac{3}{2}}\right) $.

The opposite direction, $\psi_{0,\Gamma }\left( x\right) \geq
h^{-3}\Delta _{3}^{-}\psi_{3,\Gamma }\left( x\right) $, is treated
analogously and yields the same result.

Hence,
\begin{equation*}
\psi _{0,\Gamma }\left( x\right) =\frac{x^{2}}{2}+\sum\limits_{n=1}^{M\ }%
\frac{x^{s_{n}}}{s_{n}}+O\left( x^{\frac{3}{2}}\right) \text{ \ }\left(
x\rightarrow \infty \right) \text{.}
\end{equation*}
As well-known, the latter relation implies
\begin{equation*}
\pi _{\Gamma }\left( x\right) =li\left( x^{2}\right)
+\sum\limits_{n=1}^{M}li\left( x^{s_{n}}\right) +O\left( \frac{x^{\frac{3}{2}%
}}{\log x}\right) \text{ \ }\left( x\rightarrow \infty \right) \text{.}
\end{equation*}

\textit{Proof of Corollary}. The groups $PSL\left( 2,\mathcal{O}_{K}\right) $ satisfy
the condition \eqref{e3} (see \cite[Example 4.11]{N2}). Luo, Rudnick and Sarnak \cite{LRS} proved
that the smallest eigenvalue $\lambda _{1}$ is bounded from below by $\frac{%
171}{196}$. It was further increased to $\frac{160}{169}$ by Koyama \cite{K2}.
This implies $\sum\limits_{n=1}^{M\ }\frac{x^{s_{n}}}{s_{n}}=O\left( x^{%
\frac{16}{13}}\right) $ what is obviously less than $x^{\frac{3}{2}}$. Thus,
the prime geodesic theorem takes the form
\begin{equation*}
\pi _{\Gamma }\left( x\right) =li\left( x^{2}\right) +O\left( \frac{x^{\frac{%
3}{2}}}{\log x}\right) \text{ \ }\left( x\rightarrow \infty \right) \text{.}
\end{equation*}

\section{Proof of Theorem \ref{th2}}

To reduce the exponent in the error term of the prime geodesic theorem, we
need a better control of the growth of $\sum\limits_{\left\vert
t_{n}\right\vert \leq T}\frac{x^{it_{n}}}{1+it_{n}}$. Gallagherian way to
progress in this direction is related to exclusion of a set of finite
logarithmic measure, as initially demonstrated in the classical Riemann zeta
setting \cite{G1}.

In contrast to the situation of the previous Section, the appropriate
starting point here is the explicit formula for $\psi _{2,\Gamma }$ with an
error term. A scrutiny of the argumentation in \cite{N2} brings us to the following expression
\begin{gather*}
\psi _{2,\Gamma }\left( x\right) = A_{0}x^{2}+B_{0}x^{2}\log
x+A_{1}x+A_{2}\\
+\sum\limits_{n=0}^{M}\frac{x^{s_{n}+2}}{s_{n}\left(
s_{n}+1\right) \left( s_{n}+2\right) } +\sum\limits_{n=0}^{M}\frac{x^{\widetilde{s}_{n}+2}}{\widetilde{s}_{n}\left(
\widetilde{s}_{n}+1\right) \left( \widetilde{s}_{n}+2\right) }  \\
+\sum\limits_{0\leq t_{n}\leq T}\frac{x^{s_{n}+2}}{s_{n}\left(
s_{n}+1\right) \left( s_{n}+2\right) }+\sum\limits_{0\leq t_{n}\leq T}\frac{x^{%
\widetilde{s}_{n}+2}}{\widetilde{s}_{n}\left( \widetilde{s}_{n}+1\right)
\left( \widetilde{s}_{n}+2\right) }+O\left( \frac{x^{5}}{T}\right)  \\
+\sum\limits_{\gamma _{n}\geq 0}\frac{x^{\rho _{n}+2}}{\rho _{n}\left(
\rho _{n}+1\right) \left( \rho _{n}+2\right) }+\sum\limits_{\gamma _{n}\geq 0}%
\frac{x^{\widetilde{\rho }_{n}+2}}{\widetilde{\rho }_{n}\left( \widetilde{%
\rho }_{n}+1\right) \left( \widetilde{\rho }_{n}+2\right) }\text{ \ }\left(
x\rightarrow +\infty \right) \text{.}
\end{gather*}

Let $n=\left\lfloor \log x\right\rfloor $ and take $\varepsilon >0$
arbitrarily small. For $0<Y<T$, denote by $E_{n}$ the set
\begin{center}
$\left\{ x\in %
\left[ e^{n},e^{n+1}\right) :\left\vert \sum\limits_{Y<\left\vert
t_{n}\right\vert \leq T}\frac{x^{s_{n}+2}}{s_{n}\left( s_{n}+1\right) \left(
s_{n}+2\right) }\right\vert >x^{\alpha }\left( \log x\right) ^{\beta }\left(
\log \log x\right) ^{\beta +2\varepsilon }\right\}$.
\end{center}
Then
\begin{equation}\label{e4}
\begin{gathered}
\mu ^{\times }E_{n}=\int\limits_{E_{n}}\frac{dx}{x} \\
\leq \int\limits_{e^{n}}^{e^{n+1}}\left\vert \sum\limits_{Y<\left\vert
t_{n}\right\vert \leq T}\frac{x^{s_{n}+2}}{s_{n}\left( s_{n}+1\right) \left(
s_{n}+2\right) }\right\vert ^{2}\frac{dx}{x^{1+2\alpha }\left( \log x\right)
^{2\beta }\left( \log \log x\right) ^{2\beta +4\varepsilon }} \\
\ll \frac{e^{n\left( 6-2\alpha \right)} }{n^{2\beta }\left( \log n\right)
^{2\beta +4\varepsilon }}\int\limits_{e^{n}}^{e^{n+1}}\left\vert
\sum\limits_{Y<\left\vert t_{n}\right\vert \leq T}\frac{x^{it_{n}}}{%
s_{n}\left( s_{n}+1\right) \left( s_{n}+2\right) }\right\vert ^{2}\frac{dx}{x%
}\text{.}
\end{gathered}
\end{equation}

By the Gallagher lemma \cite{G}, the last integral is dominated by
\begin{equation}\label{e5}
\int\limits_{-\infty }^{+\infty }\left( \sum\limits_{\substack{ t<\left\vert t_{n}\right\vert<t+1
\\ Y<\left\vert t_{n}\right\vert \leq T}}\frac{1}{\left\vert
s_{n}\right\vert ^{3}}\right) ^{2}dt\text{.}
\end{equation}
To proceed further, we need a more precise form of the Weyl law. If $\Gamma $
is cocompact, then $N_{\Gamma }\left( T\right) =\frac{vol\left(\Gamma \diagdown \mathbb{H}^{3}\right) }{%
6\pi ^{2}}T^{3}+O\left( \frac{T^{2}}{\log T}\right) $ \cite{DKV}. However, if $\Gamma $ is a noncompact congruence
subgroup for some imaginary quadratic number field, we still have $N_{\Gamma
}\left( T\right) =\frac{vol\left( \Gamma \diagdown \mathbb{H}^{3}\right) }{%
6\pi ^{2}}T^{3}+O\left( T^{2}\right) $ according to \cite{P}. Thus, $%
N_{\Gamma }\left( t+1\right) -N_{\Gamma }\left( t\right) =O\left(
t^{2}\right) $ in both cases. This yields
\begin{equation*}
\left( \sum\limits_{t<\left\vert t_{n}\right\vert <t+1}\frac{1}{\left\vert s_{n}\right\vert ^{3}}%
\right) ^{2}=O\left( \frac{1}{t^{2}}\right) \text{.}
\end{equation*}
Hence, the integral \eqref{e5} is bounded by $O\left( \frac{1}{Y}\right)
$. Looking back at \eqref{e4}, we obtain
\begin{equation*}
\mu ^{\times }E_{n}\ll \frac{e^{n\left( 6-2\alpha \right) }}{Yn^{2\beta
}\left( \log n\right) ^{2\beta +4\varepsilon }}\text{.}
\end{equation*}
Putting $Y=e^{n\left( 6-2\alpha \right) }n^{1-2\beta }\left( \log n\right)
^{1-2\beta +\frac{\varepsilon }{2}}$, we get
\begin{equation*}
\mu ^{\times }E_{n}\ll \frac{1}{n\left( \log n\right) ^{1+\varepsilon}}\text{.}
\end{equation*}%
Therefore, the set $E=\cup E_{n}$ has a finite logarithmic measure.

Following the same line of argumentation as in Section \ref{s3}, we derive
\begin{gather*}
\psi _{0,\Gamma }\left( x\right) \leq h^{-2}\Delta _{2}^{+}\psi
_{2,\Gamma }\left( x\right) =\sum\limits_{n=0}^{M}\frac{x^{s_{n}}}{s_{n}}%
+O\left( h^{2}\right) \\
+\frac{1}{h^{2}}\left\vert \Delta
_{2}^{+}\sum\limits_{\left\vert t_{n}\right\vert \leq T}\frac{x^{s_{n}+2}}{%
s_{n}\left( s_{n}+1\right) \left( s_{n}+2\right) }\right\vert +O\left( \frac{x^{5}}{h^{2}T}\right) \\ +\frac{1}{h}O\left( \sum\limits_{\gamma _{n}\geq 0}%
\frac{x^{\beta _{n}+1}}{\left\vert \rho _{n}\right\vert \left\vert \rho
_{n}+1\right\vert }\right) \text{,}
\end{gather*}%
where%
\begin{equation*}
\Delta _{2}^{\pm }f\left( x\right) =f\left( x\pm 2h\right) -2f\left( x\pm
h\right) +f\left( x\right) \text{.}
\end{equation*}

Now, splitting $\Delta _{2}^{+}\sum\limits_{\left\vert
t_{n}\right\vert \leq T}\frac{x^{s_{n}+2}}{s_{n}\left( s_{n}+1\right) \left(
s_{n}+2\right) }$ into
$$\Delta _{2}^{+}\sum\limits_{\left\vert
t_{n}\right\vert \leq Y}\frac{x^{s_{n}+2}}{s_{n}\left( s_{n}+1\right) \left(
s_{n}+2\right) }+\Delta _{2}^{+}\sum\limits_{Y<\left\vert t_{n}\right\vert
\leq T}\frac{x^{s_{n}+2}}{s_{n}\left( s_{n}+1\right) \left( s_{n}+2\right) }\text{,}$$
we arrive at
\begin{equation}\label{e6}
\begin{gathered}
\psi _{0,\Gamma }\left( x\right) -\sum\limits_{n=0}^{M}\frac{x^{s_{n}}}{%
s_{n}}\\
 \ll h^{2}+\sum\limits_{\left\vert s_{n}\right\vert \leq Y}\frac{x}{%
\left\vert s_{n}\right\vert }+\frac{x^{\alpha }\left( \log x\right) ^{\beta
}\left( \log \log x\right) ^{\beta +2\varepsilon }}{h^{2}}+\frac{x^{5}}{%
h^{2}T}+\frac{x^{2}}{h} \\
\ll h^{2}+xY^{2}+\frac{x^{\alpha }\left( \log x\right) ^{\beta }\left(
\log \log x\right) ^{\beta +2\varepsilon }}{h^{2}}+\frac{x^{5}}{h^{2}T}+%
\frac{x^{2}}{h}
\end{gathered}
\end{equation}
on the complement of $E$. Recall that $Y\sim x^{6-2\alpha }\left( \log x\right)
^{1-2\beta }\left( \log \log x\right) ^{1-2\beta +\frac{\varepsilon }{2}}$.
Optimizing the first three summands on the right-hand side of \eqref{e6}, we
get
\begin{equation*}
x^{\frac{\alpha }{2}}\left( \log x\right) ^{\frac{\beta }{2}}\left( \log
\log x\right) ^{\frac{\beta }{2}+\varepsilon }=x^{13-4\alpha }\left( \log
x\right) ^{2-4\beta }\left( \log \log x\right) ^{2-4\beta +\varepsilon }%
\text{.}
\end{equation*}%
This yields $\alpha =\frac{26}{9}$ and $\beta =\frac{4}{9}$. We get%
\begin{equation*}
\psi _{0,\Gamma }\left( x\right) \leq \sum\limits_{n=0}^{M}\frac{x^{s_{n}}%
}{s_{n}}+O\left( x^{\frac{13}{9}}\left( \log x\right) ^{\frac{2}{9}}\left(
\log \log x\right) ^{\frac{2}{9}+\varepsilon }\right) +O\left( \frac{x^{5}}{%
h^{2}T}\right) +O\left( \frac{x^{2}}{h}\right) \text{.}
\end{equation*}%
It is easily checked that $\frac{x^{2}}{h}$ is dominated by $%
O\left( x^{\frac{13}{9}}\right) $. Finally, $T>Y$ can be chosen so that $O\left(
\frac{x^{5}}{h^{2}T}\right) $ does not effect the bound.

The same procedure applies to $\psi _{0,\Gamma }\left( x\right) \geq
h^{-2}\Delta _{2}^{-}\psi _{2,\Gamma }\left( x\right) $. Hence,%
\begin{equation*}
\psi _{0,\Gamma }\left( x\right) =\sum\limits_{n=0}^{M}\frac{x^{s_{n}}}{%
s_{n}}+O\left( x^{\frac{13}{9}}\left( \log x\right) ^{\frac{2}{9}}\left(
\log \log x\right) ^{\frac{2}{9}+\varepsilon }\right) \text{ \ }\left(
x\rightarrow \infty ,x\notin E\right) \text{,}
\end{equation*}%
what gives us%
\begin{equation*}
\pi _{\Gamma }\left( x\right) =li\left( x^{2}\right)
+\sum\limits_{n=1}^{M}li\left( x^{s_{n}}\right) +O\left( x^{\frac{13}{9}%
}\left( \log x\right) ^{-\frac{7}{9}}\left( \log \log x\right) ^{\frac{2}{9}%
+\varepsilon }\right)\text{,}
\end{equation*}
as $x\rightarrow \infty$, $x\notin E$.

\end{document}